\documentclass[11pt]{article}

\def\href#1#2{#2}

\def\emailaddress#1{}
\def\homepage#1{}
\def\address#1{}
\def\abstract#1{}
\def\jfiftyotherinformation#1{}

\begin{document}

%
%


%

\jfiftyotherinformation{ANY OTHER INFORMATION YOU CARE TO PROVIDE}

\thispagestyle{empty}
\title{Tolstoy's Mathematics in ``War and Peace''}

\author{Paul Vit\'anyi\thanks{CWI and University of Amsterdam.
Address: CWI, Kruislaan 413, 1098 SJ Amsterdam, The Netherlands;
Email: paulv@cwi.nl; WWW:
http://www.cwi.nl/$\sim$paulv/
}}

\date{}

\maketitle


\begin{abstract}
TThe nineteenth century Russian author Leo Tolstoy
based his egalitarian views on sociology and history on mathematical
and probabilistic views, 
and he also proposed a mathematical theory of waging war.
\end{abstract}

\section{Introduction}
It is interesting to consider the excursions of mathematicians
and scientists into prose and poetry, and conversely and
less known, the explorations of poets and novelists into mathematics.

An example of the first is Luitzen E.J. Brouwer's excursion
into literature and environmentalism~\cite{Br05}, 
an appeal {\em avant la lettre} to save the earths
natural environment from human polution. In particular he wants 
to abolish the technology that enables man's supremacy over nature
and the physics and mathematics that makes this possible.
Only pure (`intuitionistic') mathematics that by its nature is unapplied and
unapplicable for evil purposes, and which is the ultimate creation
of the noble mind, should be saved. 

In another direction,
the great Russian mathematician Andrei N. Kolmogorov
was particularly interested in the 
form and structure of the poetry by the Russian author
Pushkin \cite{Ko64}. He also remarks \cite{Ko65}: ``what real meaning
is there, for example, in asking how much information is 
contained in `War and Peace'? Is it reasonable to include
this novel in the set of `possible novels,' or even to
postulate some probability distribution for this set? 
Or, on the other hand, must we assume that the individual scenes
in this book form a random sequence with `stochastic relations'
that damp out quite rapidly over a distance of several pages?

The answer to the latter question is decidly `no'. There is a
ubiquitous general theme in `War and Peace', namely, the
idea that single individuals cannot influence in any sense the course of
history (contrary to what is assumed in common history writing), 
but that the course
of history is determined by the confluence of myriad
motions of the infinitesimally small individual human acts
of free will, much as a flock of birds wheels about in unision without
any apparent governor. Here we have individual humans as interchangeable
atoms of ideal gas that in combination determine effects on macroscopic
scales such as heat and pressure, as in the nineteenth century
statistical physics of H. von Helmholtz. It serves to justify
egalitarian
doctrine. Helmholtz is also the author of the unrelated whitticism,
so true and so unknown to politicians and managers of science:
``
Whoever in the pursuit of science, seeks 
after immediate practical utility may 
rest assured that he seeks in vain,'' \cite{He62}.

It is unknown and perhaps unlikely that 
the great Russian novelist
Count Leo Nikolayevich Tolstoy 
was aware of Helmholtz's work (or directly related work); however
he may have been exposed to Laplace's \cite{La} {\em Essay} exposing
the world as a mechanistic ensemble of moving and colliding 
particles that by their combined microscopic actions produce
macroscopic effects. This is all the more likely because of
the immense influence of Laplace's work, combined with
Tolstoy's interest in justice and believable testimony in
his role as country noble passing judgement on his people.
The {\em Essay} in fact treats in great detail matters of 
believable testimony and probability of proper justice---this
was a time when relevant matters mattered.

An issue for Tolstoy is unknowability and uncertainty: he is
not really seeking a usable model so much as a reductio ad absurdum to show the
futility of the quest for explanations of wars' outcomes. 
All in all, this is not a matter of saying that the future
is in the laps of the gods, but rather that it is deterministic
and determined precisely, but practically and possibly in
principle unknowable by humans. Much like Laplace's 
viewpoint in the {\em Essay} where a demon knowing 
the positions and velocities of all particles can perfectly
predict the future and reconstruct the past, while to
the imperfect human mind not all information can be available
in a snapshot and so it is reduced to ignorance or at best
probabilistic reasoning.

The author of `War and Peace' 
had an intense interest in
mathematical approaches to the sciences,
as appears from his proposals to
found sociology, history, and the science of
war as a mathematical discipline,
much like mathematician
 John von Neumann proposed to found the science of economy as
a mathematical discipline in \cite{vN45}.

 Tolstoy's views on the matter are set forth at great length in 
`War and Peace', by many regarded
as the greatest novel in any language. Based on, or
perhaps called upon as justification for, Tolstoy's 
egalitarian philosophy, it is set
forth passionately in long interludes littered through the 
later parts of this great novel. Recall that the book is ostensibly
about the doings and adventures of a group of aristocratic people,
and in the descriptions of great battles, at the 
time of Napoleon's invasion in the bleak reaches of great Russia.
Closer inspection reveals that one of the main themes of the tale is
the insignificance and expendability of the particular
heros---like Napoleon---in the sweep of history: the events would have
infolded in the same way irrespective of the so-called main
figures. We base our treatment on Rosemary Edmonds
1957 translation into English published in Penguin Classics\cite{To69}
(part I, 1972 printing; part II, revised 1978 printing).
I will refer to the page numbers as [WP, xx].

\section{Mathematical Sociology}
Tolstoy disagrees with the view of history that ascribes the
evolution of events to individuals:

``One might have supposed that the historians, who ascribe
the actions of the masses to the will of one man,
would have found it impossible to fit the flight of Napoleon's armies
into their theory, considering that during this period of the campaign
[in Russia] the French did all they could to bring about their
own ruin, and that not a single
movement of that rabble of men $\ldots$ betrayed a hint of rhyme or reason.
But no! Mountains of volumes have been written by historians $\ldots$
[with] accounts of Napoleon's masterly arrangements and deeply
considered plans $\dots$'' [WP, 1266]

Not only that individuals cannot be the main governors of
the making of History, but:

``It is beyond the power of the human intellect to encompass
{\em all} the causes of a phenomenon.'' $\ldots$ ``the human
intellect $\ldots$ snatches at the first comprehensible 
approximation to a cause and says: `There is the cause'.''
Tolstoy goes on [WP, 1168] to explain ``in historical events
(where the actions of men form the subject of observation)
the primeval conception of a case was the will of the gods,
succeeded later on by the will of those who stand on the historical
foreground---the heroes of history.'' 

On page [WP, 1342] Tolstoy continues to unmask common misconceptions
of traditional views of History: 

``Why did things happen thus,
and not otherwise? Because they did so happen. `{\em Chance}
created the situation; {\em genius} made use of it,' says history.
But what is {\em chance}? What is {\em genius}?
The words {\em chance} and {\em genius} do not denote anything
that actually exists, and therefore they cannot be defined.
These two words merely indicate a certain degree of
comprehension of the phenomena. I do not know
why a certain event occurs; I suppose that I cannot
know: therefore I do not try to know, and I talk 
about {\em chance}. I see a force producing
effects beyond the scope of ordinary human agencies;
I do not understand why this occurs, and I cry {\em genius}.''

Now we come to the true view of history, in the spirit of the
so successful natural sciences. The ``unreasonable effectiveness
of mathematics in science'' as phrased by E. Wigner, must be
extended {\em avant la lettre} to sociology and political
history  [WP, 977]:

``To elicit the laws of history we must leave 
aside kings, ministers, and generals,
and select for study the homogeneous, infinitesemal elements which
influence the masses. No one can say how far it is possible for a man to
advance in this way to an understanding of the laws of
history; but it is obvious that this is the only path to that
end, and that the human intellect has not, so far, applied
in this direction one-millionth of the energy which historians
have devoted to describing the deeds of various kings, generals and ministers,
and propounding reflections of their own concerning those deeds.'' 

How then is this proper view of history obtained? Tolstoy discusses the
continuity of motion that was captured in laws by dividing
continuity into units. He observes that this can be done in a wrong
way,  [WP, 974]:

``Take, for instance, the well-known sophism of the ancients
which set out to prove that Achilles would never catch up with the
tortoise that had the start on him, even though Achilles traveled ten times as
fast as the tortoise: by the time Achilles has covered
the distance that separated him from the tortoise, the tortoise has advanced
one-tenth of that distance ahead of him. While Achilles does this
tenth the tortoise gains a hudredth, and so on {\em ad infinitum}.
This problem appeared to the ancients insoluble. The absurdity of the 
finding (that Achilles can never overtake the tortoise) follows from
arbitrarily separating the motion into separate units, whereas the motion
of Achilles and the tortoise was continuous. 

By adopting smaller and smaller units of motion we 
only approximate the solution of the problem but never reach it. 
It is only by admitting infinitesimal quantities and their progression
up to a tenth, and taking the sum of that geometrical progression,
that we arrive at the solution of the problem.''

Now we come to the heart of the matter: Tolstoy's proposal of 
a differential and integral analysis of history [WP, 974--975]:

``A new branch of mathematics, having attained the art of reckoning
with infinitesimal, can now yield solutions
to other more complex problems of motion which before seemed insoluble.
This new branch of mathematics, which was unknown to the 
ancients,\footnote{Apart from Archimedes and Eudoxos [PV].}
by admitting the conception, when dealing with problems
of motion, of the infinitely small and thus conforming to the
chief condition of motion (absolute continuity),
corrects the inivitable error which human intellect cannot 
but make if it considers separate units of motion instead of continuous motion.
In the investigation of the laws of historical movement
precisely the same principle operates.

The march of humanity, springing as it does from an infinite multitude
of individual wills, is continuous. The discovery of the laws of this
continuous movement is the aim of history. But to arrive at these
laws of continuous motion resulting from the sum of all those
human volitions, human reason postulates arbitrarily, separated units. The
first proceeding of a historian is to select at random a series
of successive events and examine them apart from others,
though there is and can be no {\em beginning} to any
event, for an event flows without break in continuity from another.
The second method is to study the actions of some one man---a king
or a commander---as though their actions represented the sum of many individual
wills; whereas the sum of the individual wills never finds
expression in the activity of a single historical personage.

$\ldots$ Only by assuming an infinitesimal small unit for observation---a
differential of history (that is, the common tendencies of men)---and
arriving at the art of integration (finding the sum of the infinitesimals)
can we hope to discover the laws of history.''

\section{Mathematics of War}
The causality involved in war defies simple analysis, Tolstoy says,
but is reached by the integration of the infinitesimal individual
causes,
[WP, 1184]:

``An infinite amount of freely acting forces (and nowhere is
a man freer than during a life and death struggle)
influence the course taken by a battle, and that course
can never be known beforehand
and never coincides with the direction it would have
taken under the impulsion of any single force.

If simultaneously and variously directed forces act on a given
body, the direction which that body will take cannot
be the course of anyone of the forces individually---it
will
always follow an intermediate, as it were, shortest path, or what is presented
in mechanics by the diagonal of a parallelogram of forces.'' 

In [WP, 1223---1224] Tolstoy outlines the mathematics of war and
goes into an explicit calculation that is patently false:

``Military science says, the greater the numbers [of an army]
the greater the strength. $\ldots$ For militray science to make this
assertion is like defining energy in mechanics by reference to the mass only.
It is like saying that the momenta of moving bodies will be equal or unequal
according to the equality or inequality of their masses. But
momentum (or `quantity of motion') is the product of mass and velocity.
So in warfare the strength of an army is the product of its mass
and of something else, some unknown factor $x$.'' 

He goes on to debate what this unknown $x$ may stand for and
rejects the common explanations, especially the interpretation
of $x$ as the amount of genius of the commanding general. He
goes on to say that [WP, 1224]:

``We must accept the unknown and see it for what
it is: the more or less active desire to fight and face danger.
Only then, expressing the known historical facts by means of
equations, shall we be able to compare the relative values of the
unknown factor; only then may we hope to arrive at the unknown
itself.

If ten men, batalions or divisions, fighting fifteen men, batallions
or divisions, beat the fifteen---that is, kill or capture them all
while losing four themselves, the loss will have been four on one side and
fifteen on the other. Therefore, the four are equal to the fifteen,
and we may write $4x=15y$. In other words, $x$ is to $y$
as 15 is to 4. Though this equation does not yet
give us the absolute value of the unknown factor,
it does give us a ratio between two unknowns. And by putting
a whole variety of historical data (battles, campaigns,
periods of warfare, and so on) into the form of such equations,
a series of figures will be obtained which must involve the laws
inherent in equations and will in time reveal them.''

This argument of Tolstoy is remarkable. He compares the loss
of the conquering army with the total of the vanquished army---perhaps
on the grounds that the vanquished army is totally lost.
Testing the idea by inserting
more extreme figures, such as that an army
of 1.000.000 men beats a small army of 10 men, while the conquering
army looses one man, we obtain the equation $x=10y$.
This means that the fighting spirit of the million-men
army exceeded necessarily the fighting spirit of the minuscule
ten-men army tenfold. The problem with Tolstoy's reasoning here
is that he equates the ratio of the loss of the conquering army (irrespective of
the size of the total army) and the total of the beaten army (however
small) with the ratio of the fighting spirit of the beaten army
and that of the conquering army.  In our opinion this reasoning
is hard to defend in general as is shown by substituting extrem
numbers as above. The general drift of the argument is of course
reasonable. Note that (contrary to the intention of the author)
the variables $x$ and $y$ may contain the quality of the commanders
(much as the quality of performance of a good symphony orchestra
greatly depends on the quality of the conductor).

\section{Conclusion}

It is seldom the case that a great author deems fit to incorporate
extensive discussions about mathematical foundations of social
sciences in a major literary novel. It is much more common that
scientists strive for literary redemption. Tolstoy is one of
the rare examples of the former. In fact, he gives definite proposals
to mathematize history, sociology, and the sciences of war in line
with the rational inclination of the nineteenth century.

\appendix
\section{Infinitesimal methods in Tolstoy}
Above we have extracted a few salient parts
of ``War and Peace'' advocating the use of the
methods of the calculus to study history. We used only fragments
to keep the pace; nonetheless it may be useful to provide some
more complete quotations in an appendix.
The following fragments of ``War and Peace'' are taken from the Virginia Tech
gopher site.

\subsection{Third book, Third part, Chapter 1}
Absolute continuity of motion is not comprehensible to the human
mind. Laws of motion of any kind become comprehensible to man only
when he examines arbitrarily selected elements of that motion; but
at the same time, a large proportion of human error comes from the
arbitrary division of continuous motion into discontinuous elements.
There is a well known, so-called sophism of the ancients consisting in
this, that Achilles could never catch up with a tortoise he was
following, in spite of the fact that he traveled ten times as fast
as the tortoise. By the time Achilles has covered the distance that
separated him from the tortoise, the tortoise has covered one tenth of
that distance ahead of him: when Achilles has covered that tenth,
the tortoise has covered another one hundredth, and so on forever.
This problem seemed to the ancients insoluble. The absurd answer (that
Achilles could never overtake the tortoise) resulted from this: that
motion was arbitrarily divided into discontinuous elements, whereas
the motion both of Achilles and of the tortoise was continuous.

By adopting smaller and smaller elements of motion we only
approach a solution of the problem, but never reach it. Only when we
have admitted the conception of the infinitely small, and the
resulting geometrical progression with a common ratio of one tenth,
and have found the sum of this progression to infinity, do we reach
a solution of the problem.

A modern branch of mathematics having achieved the art of dealing
with the infinitely small can now yield solutions in other more
complex problems of motion which used to appear insoluble.

This modern branch of mathematics, unknown to the ancients, when
dealing with problems of motion admits the conception of the
infinitely small, and so conforms to the chief condition of motion
(absolute continuity) and thereby corrects the inevitable error
which the human mind cannot avoid when it deals with separate elements
of motion instead of examining continuous motion.

In seeking the laws of historical movement just the same thing
happens. The movement of humanity, arising as it does from innumerable
arbitrary human wills, is continuous.

To understand the laws of this continuous movement is the aim of
history. But to arrive at these laws, resulting from the sum of all
those human wills, man's mind postulates arbitrary and disconnected
units. The first method of history is to take an arbitrarily
selected series of continuous events and examine it apart from others,
though there is and can be no beginning to any event, for one event
always flows uninterruptedly from another.

The second method is to consider the actions of some one man- a king
or a commander- as equivalent to the sum of many individual wills;
whereas the sum of individual wills is never expressed by the activity
of a single historic personage.

Historical science in its endeavor to draw nearer to truth
continually takes smaller and smaller units for examination. But
however small the units it takes, we feel that to take any unit
disconnected from others, or to assume a beginning of any
phenomenon, or to say that the will of many men is expressed by the
actions of any one historic personage, is in itself false.

It needs no critical exertion to reduce utterly to dust any
deductions drawn from history. It is merely necessary to select some
larger or smaller unit as the subject of observation- as criticism has
every right to do, seeing that whatever unit history observes must
always be arbitrarily selected.

Only by taking infinitesimally small units for observation (the
differential of history, that is, the individual tendencies of men)
and attaining to the art of integrating them (that is, finding the sum
of these infinitesimals) can we hope to arrive at the laws of history.

The first fifteen years of the nineteenth century in Europe
present an extraordinary movement of millions of people. Men leave
their customary pursuits, hasten from one side of Europe to the other,
plunder and slaughter one another, triumph and are plunged in despair,
and for some years the whole course of life is altered and presents an
intensive movement which first increases and then slackens. What was
the cause of this movement, by what laws was it governed? asks the
mind of man.

The historians, replying to this question, lay before us the sayings
and doings of a few dozen men in a building in the city of Paris,
calling these sayings and doings "the Revolution"; then they give a
detailed biography of Napoleon and of certain people favorable or
hostile to him; tell of the influence some of these people had on
others, and say: that is why this movement took place and those are
its laws.

But the mind of man not only refuses to believe this explanation,
but plainly says that this method of explanation is fallacious,
because in it a weaker phenomenon is taken as the cause of a stronger.
The sum of human wills produced the Revolution and Napoleon, and
only the sum of those wills first tolerated and then destroyed them.

"But every time there have been conquests there have been
conquerors; every time there has been a revolution in any state
there have been great men," says history. And, indeed, human reason
replies: every time conquerors appear there have been wars, but this
does not prove that the conquerors caused the wars and that it is
possible to find the laws of a war in the personal activity of a
single man. Whenever I look at my watch and its hands point to ten,
I hear the bells of the neighboring church; but because the bells
begin to ring when the hands of the clock reach ten, I have no right
to assume that the movement of the bells is caused by the position
of the hands of the watch.

Whenever I see the movement of a locomotive I hear the whistle and
see the valves opening and wheels turning; but I have no right to
conclude that the whistling and the turning of wheels are the cause of
the movement of the engine.

The peasants say that a cold wind blows in late spring because the
oaks are budding, and really every spring cold winds do blow when
the oak is budding. But though I do not know what causes the cold
winds to blow when the oak buds unfold, I cannot agree with the
peasants that the unfolding of the oak buds is the cause of the cold
wind, for the force of the wind is beyond the influence of the buds. I
see only a coincidence of occurrences such as happens with all the
phenomena of life, and I see that however much and however carefully I
observe the hands of the watch, and the valves and wheels of the
engine, and the oak, I shall not discover the cause of the bells
ringing, the engine moving, or of the winds of spring. To that I
must entirely change my point of view and study the laws of the
movement of steam, of the bells, and of the wind. History must do
the same. And attempts in this direction have already been made.

To study the laws of history we must completely change the subject
of our observation, must leave aside kings, ministers, and generals,
and the common, infinitesimally small elements by which the masses are
moved. No one can say in how far it is possible for man to advance
in this way toward an understanding of the laws of history; but it
is evident that only along that path does the possibility of
discovering the laws of history lie, and that as yet not a millionth
part as much mental effort has been applied in this direction by
historians as has been devoted to describing the actions of various
kings, commanders, and ministers and propounding the historians' own
reflections concerning these actions.

\subsection{Fourth book, Second part, Chapter 8}
A countless number of free forces (for nowhere is man freer than
during a battle, where it is a question of life and death) influence
the course taken by the fight, and that course never can be known in
advance and never coincides with the direction of any one force.

If many simultaneously and variously directed forces act on a
given body, the direction of its motion cannot coincide with any one
of those forces, but will always be a mean- what in mechanics is
represented by the diagonal of a parallelogram of forces.

If in the descriptions given by historians, especially French
ones, we find their wars and battles carried out in accordance with
previously formed plans, the only conclusion to be drawn is that those
descriptions are false.

\subsection{Second epilogue, Chapter 11}
History examines the manifestations of man's free will in connection
with the external world in time and in dependence on cause, that is,
it defines this freedom by the laws of reason, and so history is a
science only in so far as this free will is defined by those laws.

The recognition of man's free will as something capable of
influencing historical events, that is, as not subject to laws, is the
same for history as the recognition of a free force moving the
heavenly bodies would be for astronomy.

That assumption would destroy the possibility of the existence of
laws, that is, of any science whatever. If there is even a single body
moving freely, then the laws of Kepler and Newton are negatived and no
conception of the movement of the heavenly bodies any longer exists.
If any single action is due to free will, then not a single historical
law can exist, nor any conception of historical events.

For history, lines exist of the movement of human wills, one end
of which is hidden in the unknown but at the other end of which a
consciousness of man's will in the present moves in space, time, and
dependence on cause.

The more this field of motion spreads out before our eyes, the
more evident are the laws of that movement. To discover and define
those laws is the problem of history.

From the standpoint from which the science of history now regards
its subject on the path it now follows, seeking the causes of events
in man's freewill, a scientific enunciation of those laws is
impossible, for however man's free will may be restricted, as soon
as we recognize it as a force not subject to law, the existence of law
becomes impossible.

Only by reducing this element of free will to the infinitesimal,
that is, by regarding it as an infinitely small quantity, can we
convince ourselves of the absolute inaccessibility of the causes,
and then instead of seeking causes, history will take the discovery of
laws as its problem.

The search for these laws has long been begun and the new methods of
thought which history must adopt are being worked out simultaneously
with the self-destruction toward which- ever dissecting and dissecting
the causes of phenomena- the old method of history is moving.
<br>All human sciences have traveled along that path. Arriving at
infinitesimals, mathematics, the most exact of sciences, abandons
the process of analysis and enters on the new process of the
integration of unknown, infinitely small, quantities. Abandoning the
conception of cause, mathematics seeks law, that is, the property
common to all unknown, infinitely small, elements.

In another form but along the same path of reflection the other
sciences have proceeded. When Newton enunciated the law of gravity
he did not say that the sun or the earth had a property of attraction;
he said that all bodies from the largest to the smallest have the
property of attracting one another, that is, leaving aside the
question of the cause of the movement of the bodies, he expressed
the property common to all bodies from the infinitely large to the
infinitely small. The same is done by the natural sciences: leaving
aside the question of cause, they seek for laws. History stands on the
same path. And if history has for its object the study of the movement
of the nations and of humanity and not the narration of episodes in
the lives of individuals, it too, setting aside the conception of
cause, should seek the laws common to all the inseparably
interconnected infinitesimal elements of free will.

\subsection*{Acknowledgements}
I thank Peter G\'acs and Tom Koornwinder for their comments.

\end{document}